\documentclass[12pt]{amsart}
\usepackage{amscd,amsmath,amssymb,amsfonts}
\theoremstyle{plain}
\newtheorem{thm}{Theorem}
\newtheorem{lem}[thm]{Lemma}
\newtheorem{cor}[thm]{Corollary}
\newtheorem{prop}[thm]{Proposition}
\newtheorem{propdefn}[thm]{Proposition-Definition}

\theoremstyle{definition}
\newtheorem{defn}[thm]{Definition}
\newtheorem{rmk}[thm]{Remark}
\newtheorem{rmks}[thm]{Remarks}

\numberwithin{thm}{section}
\numberwithin{equation}{section}

\newcommand{\p}{\partial}

\newcommand{\eq}[2]{\begin{equation}\label{#1}#2 \end{equation}}
\newcommand{\ml}[2]{\begin{multline}\label{#1}#2 \end{multline}}
\newcommand{\ga}[2]{\begin{gather}\label{#1}#2 \end{gather}}

\newcommand{\surj}{\twoheadrightarrow}
\newcommand{\inj}{\hookrightarrow}

\newcommand{\rank}{{\rm rank}}

\newcommand{\Spec}{{\rm Spec \,}}

\newcommand{\sC}{{\mathcal C}}
\newcommand{\sD}{{\mathcal D}}
\newcommand{\sE}{{\mathcal E}}
\newcommand{\sF}{{\mathcal F}}

\newcommand{\sH}{{\mathcal H}}
\newcommand{\sI}{{\mathcal I}}

\newcommand{\sL}{{\mathcal L}}
\newcommand{\sM}{{\mathcal M}}
\newcommand{\sN}{{\mathcal N}}
\newcommand{\sO}{{\mathcal O}}

\newcommand{\sV}{{\mathcal V}}
\newcommand{\sW}{{\mathcal W}}

\newcommand{\A}{{\mathbb A}}

\newcommand{\C}{{\mathbb C}}

\newcommand{\F}{{\mathbb F}}
\newcommand{\G}{{\mathbb G}}
\renewcommand{\H}{{\mathbb H}}

\renewcommand{\P}{{\mathbb P}}

\begin{document}

\title[Fourier Transforms and rigidity]{Local Fourier Transforms and
rigidity for $\sD$-Modules}
\author{Spencer Bloch}
\address{Dept. of Mathematics,
University of Chicago,
Chicago, IL 60637,
USA}
\email{bloch@math.uchicago.edu}

\author{H\'el\`ene Esnault}
\address{Mathematik,
Universit\"at Essen, FB6, Mathematik, 45117 Essen, Germany}
\email{esnault@uni-essen.de}
\date{December 17, 2003}
\begin{abstract}
Local Fourier transforms, analogous to the $\ell$-adic local Fourier
transforms \cite{Lau}, are constructed for connections over $k((t))$.
Following a program of Katz \cite{Ka}, a meromorphic connection on a
curve is shown to be rigid, i.e. determined by local data at the
singularities, if and only if a certain infinitesimal rigidity condition
is satisfied. As in \cite{Ka}, the argument uses local Fourier transforms
to prove an invariance result for the rigidity index under global Fourier
transform. A key technical tool is the notion of good lattice pairs for a
connection \cite{DPSR}. 
\end{abstract}
\maketitle
\begin{quote}

\end{quote}

\section{Introduction}

In an important article, G. Laumon \cite{Lau} applied the $\ell$-adic
Fourier transform to study epsilon factors associated to $\ell$-adic
sheaves on curves over finite fields. As a key tool, he
defined local Fourier transforms $\sF(0,\infty),
\sF(\infty, 0),\sF(\infty,\infty)$ for $\ell$-adic sheaves on $\Spec
\F_q((t))$. Recently, we applied his ideas to study epsilon factors
associated to holonomic $\sD$-modules on curves. The purpose of this
paper is to develop local Fourier transforms for meromorphic connections
over Laurent series fields. We show that these have properties precisely 
analogous to the Laumon local $\ell$-adic local Fourier transforms.

As an application of our construction, we consider 
the index of rigidity of meromorphic connections 
as defined by N. Katz \cite{Ka}. 
A local system 
 on $\P^1\setminus S$ over the complex numbers, 
where $S$ is a finite set of points, 
is rigid when it is uniquely determined by its local monodromies. 
Similarly, an $\ell$-adic representation over $\P^1\setminus S$ over a
finite field is rigid if it is uniquely determined by its restriction
to the Laurent power series fields at the singularities. N. Katz defines
 those concepts in \cite{Ka} and shows a fundamental classification
 theorem: rigid systems of tame $\ell$-adic representations
 are all obtained from
rank one $\ell$-adic sheaves by applying Fourier transform and
convolution (\cite{Ka}, Main  Theorem 5.2.1). A similar result holds 
for complex local systems.
One of Katz's key ideas is to define
an index of rigidity \cite{Ka} Chapter 3,  and to show that under a 
suitable 
assumption on the system or $\ell$-adic representation, this index is
invariant under Fourier transform. His proof relies on a theorem
of Laumon on the local Fourier transform and 
the equivalence of categories it yields between local $\ell$-adic
representations at $t=0$ and at $t'=\infty$ with slopes $<1$ and
between local $\ell$-adic reprenstations at $\infty$ with slopes $>1$
(\cite{Ka}, Theorem 3.0.2). Katz raises the question (op. cit, p.
10) whether an analogous invariance under Fourier is true for the index of
rigidity of a meromorphic connection on $\P^1$.  We show using the local
Fourier transforms and Katz's arguments that this is the case  (see
Theorem \ref{thm4.2}). As an application we show that an irreducible
meromorphic connection is rigid if and only if it has index of rigidity
$2$. One may hope  to use these ideas to classify rigid meromorphic
connections with irregular singular points, though this question is not
addressed here.

We remark that in \cite{Mal3}, Malgrange constructs a complex 
analytic microlocalization for analytic connections which yields a
microanalytic construction of
$\sF(0, \infty)$ and $\sF(\infty, 0)$. 

Finally we show that Katz' index of rigidity is equivalent to rigidity when the connection is irreducible. This cohomological characterization
 is not known for $\ell$-adic representations. 
\\ \ \\
{\it Acknowledgements:} We have greatly benefited from correspondence
with Nick Katz on his theorem. We also thank Alexander Beilinson for helpful 
comments and correspondence. 

\section{Grothendieck's theorem on formal cohomology}\label{sect2}

\begin{thm}[\cite{Gr}, Th\'eor\`eme 4.1.5)]\label{thm:gro}: 
Let $f: X\to S$ be a proper morphism of noetherian schemes, let 
$S'\subset S$ be a closed subset 
 defined by the ideal sheaf $\sI\subset \sO_S$,   
and let 
 $X'=f^{-1}(S')$. Then if $\sF$ is a coherent sheaf on $X$, the natural 
restriction map 
\ga{}{\varprojlim_{\ell} (R^nf_*\sF)\otimes_{\sO_S}\sO_S/\sI^\ell \to 
\varprojlim_{\ell} R^nf_*(\sF\otimes_{\sO_X} \sO_X/f^*\sI^\ell)\notag}
is an isomorphism for all $n\ge 0$. Both sides coincide with the formal
cohomology $R^n\widehat f_* \widehat \sF$ on $\widehat S$, where
$\ \widehat{}\ $ refers to the formal scheme completions along $S'$ and
$f^{-1}S'$. 
\end{thm}
We will apply this theorem to the following situation. 
Let $k$ be a field of characteristic 0, and let $(M, \nabla)$ be a 
connection on 
$\A^1\setminus S$, where $S$ is a finite collection of points. 
We set $j:\A^1\setminus S\to \A^1, k: \A^1\to \P^1$ for the open 
embeddings. 
We consider the projection $p_2: X=\P^1\times \P^1\to \P^1$ 
and to distinguish the two factors, 
we denote by $t$ the parameter on the left $\A^1$ and by $t'$ the 
parameter on the right $\A^1$. We consider the rank 1 connection
$\psi=(\sO_{\A^1\times \A^1}, d+d(tt'))$ on $\A^1\times \A^1$ and
consider the cohomology
\ga{2.1}{R^1p_{2*}\Big(p_1^* k_*j_*M\xrightarrow{p_1^*\nabla\otimes \psi}
p_1^*(\omega(*T)\otimes  k_*j_*M)\Big),}
with $T=S\cup \infty$. In computing this cohomology, $\psi$ is viewed as
a relative connection with operator $d_t + t'dt$ which is extended
meromorphically to $t=\infty$.

Write $\Omega$ for the direct image of $\Omega^1_{\A^1\times \A^1}$ on
$\P^1\times \P^1$. The standard diagram, where the middle column
calculates the de Rham cohomology on $\P^1\times \P^1$,  
(for which 
$\psi$ is viewed as a connection relative to $k$), and the right
hand column is de Rham cohomology relative to $p_2$
\begin{tiny}
\minCDarrowwidth.1cm
\ga{2.2}{
\begin{CD}
@.  p_1^* k_*j_*M @>\cong >>  p_1^* k_*j_*M\\
@.  @V p_1^*\nabla \otimes \psi VV @V p_1^*\nabla \otimes \psi VV\\
p_2^*\omega(*(0+\infty)) \otimes p_1^* k_*j_*M
@>>>\Omega(p_1^*(*T))\otimes
 p_1^* k_*j_*M @>>>  p_1^*(\omega(*T)\otimes  k_*j_*M)\\
@V 1\otimes p_1^*\nabla \otimes \psi VV @V p_1^*\nabla \otimes \psi VV \\
p_2^*\omega(*(0+\infty)) \otimes p_1^*(\omega(*T)\otimes  k_*j_*M)
@>> \cong > \makebox[2cm][l]{$p_2^*\omega(*(0+\infty)) \otimes 
p_1^*(\omega(*T)\otimes  k_*j_*M)$}
\end{CD}
}
\end{tiny}
yields, via the connecting homomorphism, a connection 
\ga{2.3}{ \eqref{2.1} \to \omega(*(0+\infty))\otimes \eqref{2.1}.}
We want to study the formal completion of 
this connection (which is the Fourier 
transform of $M$, see section \ref{sec:lft}) 
at the points $t'=\infty$ and $t'=0$.

\begin{cor}\label{cor:gro} 
Assume given on $\P^1$ vector bundles $\sV_i\subset k_*j_*M$ 
with the property that 
$(p_1^*\nabla \otimes \psi )(p_1^*\sV_1)\subset p_1^*\omega(T)\otimes
p_1^*\sV_2$ so that the inclusion of complexes
\ga{}{\Big(p_1^*\sV_1\xrightarrow{p_1^*\nabla\otimes \psi}
p_1^*\omega(T) 
\otimes p_1^*\sV_2\Big)
\subset \notag \\
\Big(p_1^* k_*j_*M\xrightarrow{p_1^*\nabla\otimes \psi}
p_1^*(\omega(*T))\otimes  k_*j_*M\Big)\notag}
is a quasi-isomorphism. 
Then one has
\ga{}{
\eqref{2.1}\otimes_{\sO_{\P^1}} k((u'))= \notag\\
\H^1\Big(\P^1[[u']], \sV_1[[u']] \xrightarrow{\nabla +t'dt} 
\omega(T)\otimes \sV_2[[u']]\Big)\otimes_{k[[u']]} k((u'))\notag  
}
where $u'=t'-a$ for some point $a$ or $u'=\frac{1}{t'}$. Here
$\P^1[[u']]$ is the formal scheme obtained by completing $\P^1\times_k
\Spec k[[u']]$ at the central fibre.
\end{cor}
\begin{proof}
Indeed, the $E_1$  spectral sequence
$E^{ab}=R^b(p_2)_*K^a \Rightarrow R^{a+b}(p_2)_* K^\bullet$ for the complex
$K^\bullet: p_1^*\sV_1 \to p_1^*\omega(T)\otimes p_1^*\sV_2$ yields
\ga{2.4}{\eqref{2.1}\otimes_{\sO_{\P^1}} k[[u']]=
\H^1\Big(\P^1[[u']], \sV_1[[u']] \xrightarrow{\nabla +t'dt} 
\omega(T)\otimes \sV_2[[u']]\Big).}
The assertion of the corollary follows by tensoring with
$\otimes_{k[[u']]} k((u'))$.  
\end{proof}
\begin{cor} \label{cor:gro2}
Let $K=k(t')$. With the assumptions as in Corollary \ref{cor:gro}, 
one has 
\ga{}{\H^1\Big(\P^1\times_k K, p_1^*k_*j_*M \xrightarrow{p_1^*\nabla +
\psi} p_1^*\omega(T)\otimes k_*j_*M\Big)\otimes_K k((u'))= \notag \\
\H^1\Big(\P^1[[u']], \sV_1[[u']] \xrightarrow{\nabla +t'dt} 
\omega(T)\otimes \sV_2[[u']]\Big)\otimes_{k[[u']]} k((u')) \notag}
and the latter does not depend on the choice of $\sV_i$ as in Corollary
 \ref{cor:gro}.
\end{cor}
The aim of section \ref{sec:lft} 
will be in particular to show the existence of such $\sV_i$. 

In general, for a bounded below complex of sheaves $\sC$ on a
topological space $X$, one has a spectral sequence
$E_2^{p,q}=H^q(X,\sH^p) \Rightarrow \H^{p+q}(X,\sC)$ where $\sH^p$ is
the $p$-th cohomology sheaf of $\sC$.  We can apply this with $X$ the
formal scheme $\P^1[[z']]$ with $z'=\frac{1}{t'}$, and $\sC$ the complex 
\eq{2.5}{\sV_1[[z']] \xrightarrow{z'\nabla +dt} \omega(T)\otimes
\sV_2[[z']]
} 
as above, placed in degrees $[0,1]$. In this case, the
differential is easily seen to be injective, so the hypercohomology in
degree $1$ is given by $H^0(\P^1[[z']], \sH^1)$. 
\begin{lem} The sheaf $\sH^1$ in this case is supported at the points of
$\P^1_k = \P^1[[z']]$ where $\nabla$ has singularities. At a point $s \in
S \subset \A^1$ where $M$ has irregularity $n$, $\sH^1_s$ is a free
$k[[z']]$-module of rank $\rank M + n$. If $M$ is smooth or has a regular
singular point at $\infty$, then $\sH^1_\infty = (0)$. 
\end{lem}
\begin{proof} For $x\in \A^1$ any point at finite distance, choose a
$k$-vector space complement $U$
\eq{2.6}{(\omega(T)\otimes\sV_2)_x = (\sV_{1,x}\wedge dt) \oplus U.
}
It is straightforward to identify $\sH^1_x \cong U[[z']]$. The assertions
for $x$ now follow from Deligne's theory of good lattices (\cite{DPSR},
lemme 6.21).  In particular, $U$ is a finite dimensional $k$-vector space
of dimension
${\rm dim} \frac{1}{u}\sV_2/\sV_2 +{\rm dim}\sV_2/\sV_1= {\rm rank}(M) +n$ 
(see the
discussion in section \ref{sec:lft}, particularly formula 
\eqref{3.2} and Proposition \ref{prop3.13}, (i).) When $M$
has at worst a regular singular point at $\infty$, one can take
$\sV_{2,\infty} =
\sV_{1,\infty}(\infty)$. Since $\infty \in T$, in this case
$\sV_{1,\infty}\wedge dt = (\omega(T)\otimes\sV_2)_\infty$, and one
concludes by a variant of the above argument. 
\end{proof}

Let $x\in \A^1$ be as above, and write $\widehat \sV_{i,x} = \sV_{i,x}
\otimes \widehat \sO_x$ for the formal completion. We have
\eq{}{(\omega(T)\otimes\sV_2)_x\Big/\sV_{1,x}\wedge dt \cong
(\omega(T)\otimes\widehat \sV_2)_x\Big/\widehat\sV_{1,x}\wedge dt. 
}
We conclude
\begin{cor}\label{cor2.5} With notation as above
\eq{}{\sH^1_x \cong \rm{coker}\Big(\widehat\sV_{1,x}[[z']]
\xrightarrow{z'\nabla +dt} 
\omega(T)\otimes \widehat\sV_{2,x}[[z']] \Big).
}
Moreover, $\sH^1_x\otimes_{k[[z']]} k((z'))$ depends only on the formal
meromorphic connection $M\otimes k((t_x))$, where $t_x$ is a local
parameter at $x$.
\end{cor}
\begin{proof} The last assertion follows from Corollary \ref{cor:gro2}
(with $u'=z'$) and the spectral sequence. 
\end{proof}
\begin{rmk} \label{rmk2.6}
Let $x\neq 0$. Let us consider the Fourier transform of $M$ ``centered at'' $x$, that is consider the definition \eqref{2.5} with $t$ replaced by $t_x$. Call
$\sH^1(t_x)$ the hypercohomology sheaf. Then one obviously has 
\ga{2.10}{\sH^1(t_x)_0\otimes xd(\frac{1}{z'})=
\sH^1_x}
where $xd(\frac{1}{z'})$ is the connection on $k((z'))$ which to 1 assigns 
$xd(\frac{1}{z'})$.
\end{rmk}

\section{Local Fourier transforms} \label{sec:lft}

In \cite{Lau}, section 2, G. Laumon defines the local 
Fourier transforms of an $\ell$-adic representation over
$\Spec \F_q((t))$.  The aim of this section 
is to  define the corresponding notion for connections on 
$\Spec k((t))$, where $k$ is a field of characteristic 0. 

Let $k$ be a field of characteristic zero, and let $\sM$ be a holonomic 
$\sD=k[t,\p_t]$-module over the affine line $\A^1={\rm Spec}k[t]$. We
recall the definition (see \cite{KaLa}, or \cite{Mal} chapter VI) 
\begin{defn} \label{defn3.1} The Fourier transform $\sF(\sM)$ of $\sM$
is the $\sD'=k[t', \p_{t'}]$-module obtained
by keeping the same $k$-vector space $\sM$ but setting 
\begin{itemize}
\item[i)]
$t'\cdot m= -\p_t \cdot m$
\item[ii)] 
$\p_{t'}\cdot m= t\cdot m$ for all $m\in \sM$. 
\end{itemize}
\end{defn}

\begin{lem} \label{lem3.2} $\sF(\sM)$ is the Gau{\ss}-Manin connection 
on $$H^1(\sM\otimes_k k[t']\xrightarrow{\p_t + t'} \sM\otimes_k k[t']).$$ 
\end{lem} 
\begin{proof}
One has the diagram
\ga{3.1}{\begin{CD}
@. @. \sM @>=>> \sM \\
@. @. @V \alpha VV @V \beta VV \\
0 @>>> \sM\otimes_k k[t'] @>\p_t+t'>> \sM\otimes_k k[t']@>p>> H^1 @>>> 0
\end{CD}
}
where $p$ is the quotient map to the first cohomology of $\p_t+t'$, 
$\alpha(m)=m\otimes 1, \ \beta= p\circ \alpha$. We first note that $\beta$ is an isomorphism of $k$-vector spaces. 
Indeed, $\beta$ is injective as $m\otimes 1$ can't be in the image of 
$\p_t+t'$. Given 
$\mu=\sum_{i=0}^N m_i (t')^i, m_N\neq 0, N\ge 1$, 
then $ \mu - (\p_t+t')(m_N (t')^{N-1})$ has degree $\le (N-1)$ in $t'$. 
Inductively, one sees that every class in $H^1$ is the class of
some  $m\otimes 1$ so $\beta$ is onto. Now $p(m\otimes t')=-p(\p_tm\otimes
1)$ which shows the  relation  i) of Definition \ref{defn3.1}.
To see ii), one computes the Gau{\ss}-Manin connection on $p(m\otimes 1)$. 
One has $((\nabla+ d(tt'))/dt)(m\otimes 1)=(\p_t+t')(m\otimes 1) + 
(tm\otimes 1)dt') \equiv (tm\otimes 1)dt'$, implying ii). 
\end{proof}

To calculate Fourier transforms we will use Deligne's {\it good}
lattices for irregular connections, defined in \cite{DPSR}, Lemme 6.21.
Let $X/k$ be a smooth curve over, $j: X\subset \bar{X},
\Sigma:=\bar{X}\setminus X$ be a smooth compactification, and let $\sM$
be a smooth connection on $X$, meromorphic along $\Sigma$. By definition,
a {\it pair of good lattices} $\sV, \sW\subset j_*\sM$  is a pair of vector
bundles on $\bar{X}$ satisfying the following conditions
\begin{itemize}
\item[1)] $\sV\subset \sW\subset j_*\sM$
\item[2)] $\nabla(\sV)\subset \omega_{\bar{X}}(\Sigma)\otimes \sW$
\item[3)] the inclusion of complexes
$$(\sV\xrightarrow{\nabla} \omega_{\bar{X}}(\Sigma)\otimes \sW)\to
(j_*\sM\xrightarrow{\nabla} \omega_{\bar{X}}\otimes j_*\sM)$$
is a quasi-isomorphism.
\end{itemize}

Notice that these conditions are purely local. For $\sigma \in \Sigma$,
let $t_\sigma$ be a local parameter at $\sigma$. It suffices to construct
lattices for $\sM \otimes k((t_\sigma))$ satisfying the analogous
conditions.

Deligne (\cite{DPSR}, p.110--112) shows the existence of good lattices. If $\sV, \sW$ are good
lattices, so are $\sV(D), \sW(D)$ for any divisor $D$ supported on
$\Sigma$.  At a point $\sigma \in \Sigma$  which is regular singular, one
has $\sV\otimes
\sO_{\bar{X}, \sigma} =
\sW\otimes \sO_{\bar{X}, \sigma}$.

For $\sigma \in \Sigma$, the dimension of the finite dimensional vector
space 
\eq{3.2}{\sW
\otimes \sO_{\bar{X}, \sigma}/\sV\otimes \sO_{\bar{X}, \sigma}
} is 
independent of the choice of $\sV, \sW$ and is equal to the 
{\it irregularity} of $\sM$ at $\sigma$. 

\begin{lem}\label{lem3.3} Let $M$ be a connection on
$k((t))$. Then  the slopes of $M$ are $\le 1$ (resp. $\ge 1$) if and only
if 
 there
exists a pair $\sV, \sW$ of good lattices such that $\sV \subset 
\sW\subset \sV(0) $ (resp. $\sV \subset \sV(0) \subset \sW$). Moreover for such a pair of good lattices, $\sV(0)=\sW$ if and only if the slopes are $=1$. 
\end{lem}
\begin{proof} We may assume $M$ is indecomposible, i.e. not of the form
$M_1 \oplus M_2$ for $M_i \neq 0$. It follows (\cite{Mal}, Th. 1.5, p. 45)
that $M$ has a single slope. The if part is clear. 
We prove necessity. If the slope is zero, the connection is
regular singular. Then one has Deligne's lattices $\sV$ (\cite{DPSR}, 
Th\'eor\`eme 4.1 and Corollaire 3.14) with respect to which the
connection has logarithmic poles thus $\sV=\sW$.  Assume the slope is
$>0$. As in \cite{BBE}, section 5.9, we may assume
$M = N\otimes U$, where $U$ is regular singular with unipotent monodromy
and $N = \pi_*L$ where $L$ is a rank $1$ connection on a finite covering
$\pi: \Spec K \to \Spec k((t))$. The integral closure of $k[[t]]$ in $K$
is a complete, equicharacteristic $0$ discrete valuation ring, so it has
the form $k'[[u]]$ for some $k'/k$ finite (\cite{Z}, Cor. 2, p. 280).
Here $u$ satisfies an Eisenstein polynomial of some degree $p$ over
$k((t))$. If we fix a trivialization $L= k'((u))\cdot e$ and write
$\nabla(e) = e\otimes (a_{-n}u^{-n} + a_{-n+1}u^{-n+1} +
\ldots)\frac{du}{u}$ with $n\ge 1$, then we get a good lattice pair for
$L$ taking $\sV_L = k'[[u]]\cdot e$ and $\sW_L = u^{-n}\sV_L$. It is now
straightforward to check that $(\pi_*\sV_L \otimes U, \pi_*\sW_L \otimes
U)$ is a good lattice pair for $M$, and that $M$ has slope $\frac{n}{p}$.
When $\frac{n}{p}\le 1$ (resp. $\frac{n}{p}\ge 1$) we have $\sW_L =
u^{-n}\sV_L \subset t^{-1}\sV_L$ (resp. $\sW_L =
u^{-n}\sV_L \supset t^{-1}\sV_L$). The assertion of the lemma follows by
applying $\pi_*$ and tensoring with $U$. 
\end{proof}

\begin{rmk}
It is not true that if the slope condition is as in Lemma \ref{lem3.3}, then all pairs of good lattices fulfill the relations of the lemma. Indeed, the ones constructed by Deligne \cite{DPSR}, Lemme 6.21 do not always. 
For example, the connection on $\oplus_1^3 \sO$ with connection matrix
\ga{}{\left(\begin{matrix}
0 & 0 & \frac{1}{t^2}\\
1 & 0 & 0 \\
0 & 1 & 0 \end{matrix} \right)\frac{dt}{t}}
does not. 
The slope is $\frac{2}{3}<1$, $\sV$ as a $\sO$-module is generated by 
$e_1, 
e_2, e_3$ while $\sW$ is generated by $e_2, e_3, \frac{1}{t^2}e_1$. 
\end{rmk}

Let  $M$ be a connection on $k((t))$. Recall from 
\cite{KaD} Theorem (2.4.10) that there is a canonical (but not
unique) functorial smooth extention
$\sM$ to $\G_m = \Spec k[t,t^{-1}]$ which has regular singular points at
$t=\infty$. We shall refer to $\sM$ as the {\it Katz extension} of $M$. Our
arguments will use the existence of a Katz extension, but nothing about
its properties. 

\begin{propdefn}[Local Fourier from $0$ to $\infty$] \label{propdefn3.5}
Let $M$ be a connection on $k((t))$ and let $\sM$ be the Katz extension
to a meromorphic connection on $\P^1$ with regular singular point at
$\ t=\infty$. Let $t'$ be the Fourier transform coordinate, and write
$z'=\frac{1}{t'}$. Then the Fourier transform connection (Definition
\ref{defn3.1}) restricted to the Laurent series field at $t'=\infty$,
$\sF(\sM)\otimes_{k[t']} k((z'))$ on
$k((z'))$ depends only on
$M$ and not on the choice of $\sM$. We call it the local Fourier transform 
of $M$ from 0 to $\infty$ and denote it by $\sF(0,\infty)(M)$.
Concretely, if $\widehat \sV,\ \widehat\sW$ is a good lattice pair for
the formal connection $M$, 
\eq{3.4}{\sF(0,\infty)(M) = \rm{coker}\Big(\widehat \sV((z'))
\xrightarrow{z'\partial_t + 1} \frac{1}{t}\widehat \sW((z'))\Big).
}
\end{propdefn}
\begin{proof}We apply the discussion of section \ref{sect2} to the
Katz extension $\sM$ of $M$. Thus, with notation as in
Corollary \ref{cor2.5}, $\sF(0,\infty)(M) := \sH^1_0 \otimes k((z'))$.
Independence of choice of good lattices follows from that corollary. 
\end{proof}

\begin{rmk} If the connection $M$ on $\Spec k((t))$ extends smoothly
across $\Spec k[[t]]$, $\sF(\sM)$ is supported at $t'=0$ and
$\sF(0,\infty)(M) = (0)$. For this reason, we will assume when working
with the local Fourier transform that $M^\nabla = (0)$. 
\end{rmk}

The construction of $\sF(0,\infty)(M)$ is independent of the choice of a
good lattice pair. In particular, we can take $\widehat \sW$ as large as
we like. The composition 
\eq{3.5}{\frac{1}{t}\widehat\sW \inj \frac{1}{t}\widehat\sW((z')) \surj
\sF(0,\infty)(M)  
}
therefore extends to a $k$-linear map
\eq{3.6}{\iota: M \to \sF(0,\infty)(M).
}
\begin{prop}\label{prop3.7} Assume $M^\nabla = (0)$. Then $\iota$ is an
isomorphism of
$k$-vector spaces. One has
$\iota \circ \p_t= -\frac{1}{z'}\circ \iota $ and 
$\iota \circ t = -(z')^2\p_{z'} \circ \iota$. 
\end{prop}
\begin{proof} Suppose $\iota(m) = 0$. Taking $\widehat\sW$ to be large,
we can assume $m\in \frac{1}{t}\widehat\sW$. 
Then $m$ has to be of the shape $(z'\p_t+1)(\sum_{\ell = N}^\infty 
(z')^\ell v_\ell)= (z')^N v_N +\sum_{\ell \ge N+1} (z')^{\ell} (\p_t v_{\ell-1} +
v_\ell) $ for some $v_\ell \in \widehat\sV$. This implies $N=0,   (\p_t)^\ell
m=(-1)^\ell v_\ell \in \widehat\sV$ for $\ell \ge 0$ so
\eq{}{m = (z'\partial_t+1)\Big(\sum_{\ell=0}^\infty (-1)^\ell z'{}^\ell
\partial_t^\ell
m \Big).
}
But the fact that $(\p_t)^\ell m \in \widehat\sV$ means the
$k[[t]][\partial_t]$-submodule of $M$ generated by
$m$ is finitely generated as a $k[[t]]$-module and so necessarily has a
horizontal section, contradicting our assumption. 

The assertion that $\iota\circ\partial_t = -\frac{1}{z'}\circ\iota$ is
clear. 
Our assumption that $M^\nabla = (0)$ implies $\partial_t : M \cong
M$  (cf. for example \cite{Mal2}, Thm. 2.1 (b)), so we may speak of
$\partial_t^{-1}$. This operator is $t$-adically contracting. Indeed,
$\widehat\sV \subset\widehat\sW$ and $\partial_t : \widehat\sV
\cong \frac{1}{t}\widehat\sW$. Clearly $\iota\circ\partial_t^{-1} =
-z'\circ\iota$. It follows that the image of $\iota$ is closed under
taking Laurent series in $z'$, from which surjectivity is clear.
As for $\iota \circ t=-(z')^2\p_{z'}\circ \iota$, 
the computation is as in 
Lemma \ref{lem3.2}. 
\end{proof}

\begin{defn}[Local Fourier from $\infty$ to $0$]\label{defn3.8} Let $M$
be a connection on $k((z))$. Assume $M^\nabla = (0)$ and that the slopes
(\cite{Mal}, chap. III) are all $<1$. Let $\widehat\sV, \widehat\sW$ be a
pair of good lattices for $M$. By Lemma \ref{lem3.3} we may assume
$z^2\partial_z\widehat\sV \subset \widehat\sV$. Then
\eq{}{\sF(\infty,0)(M) := \text{coker}\Big(\widehat\sV((t'))
\xrightarrow{-z^2\partial_z + t'} \widehat\sV((t'))\Big). 
}
\end{defn}

\begin{prop}\label{prop3.9}Let $M$ be as in the definition. Then
$\sF(\infty,0)(M)$ is independent of the choice of good lattices. The
natural map
\eq{}{\iota: \widehat\sV \inj \widehat\sV((t')) \surj \sF(\infty,0)(M)
}
extends to an isomorphism of $k$-vector spaces $\iota: M \cong
\sF(\infty,0)(M)$. We have $\iota\circ z^2\partial_z = t'\circ\iota$ and
$\iota\circ \frac{1}{z}= -\partial_{t'}\circ\iota$
\end{prop}
\begin{proof}To show independence of the choice of lattices, let
$\widehat\sV, \widehat\sW$ and $\widehat\sV',\widehat\sW'$ be two pairs
of good lattices with $\widehat\sV,
\widehat\sV'$ stable under $z^2\partial_z$. We may assume $\widehat\sV
\subset \widehat\sV'$, and we have to show $-z^2\partial_z + t'$ is
an isomorphism 
on $(\widehat\sV'/\widehat\sV)((t'))$. But this is clear
because the slope condition forces $z^2\partial_z$ to be nilpotent. 

Because $z^2\partial_z$ is injective on $M$, an equation of the form 
\eq{}{m = (-z^2\partial_z + t')\sum_{n=N}^\infty  v_nt'{}^n
}
forces $m = -(z^2\partial_z )^{n+1}v_n, \ v_n \in \widehat\sV$. Again by
the slope condition, this forces $m=0$, so $M \inj \sF(\infty,0)(M)$. 

The identity $\iota\circ z^2\partial_z = t'\circ\iota$ is clear from the
definition. Since $z^2\partial_z$ is bijective and contracting, it
follows that the image of $\iota$ is closed under taking Laurent series
in $t'$, so $\iota$ is surjective as well. 
 The proof that $\iota \circ \frac{1}{z}=\p_{t'}\circ \iota$, is as
in  Lemma \ref{lem3.2}. 
\end{proof}

\begin{prop}\label{prop3.10} The functors $\sF(0,\infty)$ and
$\sF(\infty,0)$ are inverses and define equivalences of categories
\ml{}{\Big\{k((t))-\text{connections with no horiz. sects.}\Big\}
\leftrightarrow \\
\Big\{k((z'))-\text{connections with no
horiz. sects. and slopes $<1$}\Big\}  }
\end{prop}
\begin{proof} Let $M$ be a $k((t))$-connection with $M^\nabla =
(0)$. By Proposition \ref{prop3.7},  $M \cong \sF(0,\infty)(M)$ and the
action of $-z'{}^2\partial_{z'}$ on $\sF(0,\infty)(M)$ corresponds to
multiplication by $t$ on $M$. In particular, $\sF(0,\infty)(M)^\nabla =
(0)$. Similarly, given $N$ a $k((z'))$-connection with no horizontal
sections and slopes $<1$, we have by Proposition \ref{prop3.9}, $N \cong
\sF(\infty,0)(N)$ and $\partial_{t}$ on $\sF(\infty,0)(N)$ corresponds
to multiplication by $-\frac{1}{z'}$ on $N$, so $\sF(\infty,0)(N)$ has no
global sections and the functors are defined. Finally, under the vector
space identifications, the operators intertwine as indicated: 
\ga{}{M \cong \sF(0,\infty)(M) \cong
\sF(\infty,0)\Big(\sF(0,\infty)(M)\Big) \\
t\qquad\qquad - z^2\partial_z \qquad\qquad\qquad - t\qquad\qquad\qquad \notag
\\
\partial_t \qquad\qquad -\frac{1}{z} \qquad\quad\qquad\quad -\partial_t
\qquad\qquad\qquad\notag 
}
It follows that $\sF(\infty,0)\circ\sF(0,\infty) = [t\mapsto -t]^*$. The
argument in the other direction is similar. 
\end{proof}

\begin{defn} \label{defn3.11} Let $M$ be a connection on $k((z))$ and assume all slopes of
$M$ are $>1$. Let $\widehat\sV, \widehat\sW$ be a good lattice pair with
$\frac{1}{z}\widehat\sV \subset \widehat\sW$. Define
\eq{3.13a}{\sF(\infty,\infty)(M) := \text{coker}\Big(\widehat\sV((z'))
\xrightarrow{-z'z^2\partial_z+1} z\widehat\sW((z'))\Big).
}
$\sF(\infty,\infty)(M)$ is a connection on $k((z'))$. 
\end{defn}

\begin{prop}\label{prop3.12} Let $M$ be a connection on $k((z))$ with
slopes
$>1$. 
\begin{enumerate}
\item[(i)] $\sF(\infty,\infty)(M)$ is independent of the choice of good
lattice pair.
\item[(ii)] The evident projection $z\widehat\sW \to
\sF(\infty,\infty)(M)$ extends to an isomorphism of $k$-vector spaces
$\iota : M \cong \sF(\infty,\infty)(M)$.
\item[(iii)] The operators $z^2\partial_z$ and $-\frac{1}{z}$ on $M$ coincide
with the operators $\frac{1}{z'}$ and $z'{}^2\partial_{z'}$ on
$\sF(\infty,\infty)(M)$. In particular, $\sF(\infty,\infty)(M)$ has all
slopes $>1$. 
\item[(iv)] $\sF(\infty,\infty)\circ\sF(\infty,\infty)(M) \cong
\sigma^*M$, where $\sigma: k((z)) \to k((z))$ is the automorphism
$z\mapsto -z$. In particular, $\sF(\infty, \infty)$ is an
auto-equivalence of the category of connections on $k((z))$ with
slopes $>1$. 
\item[(v)] Let $\sM$ be a Katz extension of $M$ to a meromorphic
connection on $\P^1$, smooth over $\G_m$ with a regular singular point at
$z=\infty$. Then $\sF(\infty,\infty)(M) \cong \sF(\sM)\otimes k((z'))$.
\end{enumerate}
\end{prop}
\begin{proof}\begin{enumerate}
\item[(i)] It suffices to consider good lattice pairs $\widehat\sV,
\widehat\sW$ and $\widehat\sV', \widehat\sW'$ with $\widehat\sV \subset
\widehat\sV'$ and $\widehat\sW \subset\widehat\sW'$. Again by
\cite{Mal2}, Thm. 2.1 (b), $z\partial_z : \widehat\sV \cong
\widehat\sW$ (resp. $\widehat\sV' \cong \widehat\sW'$). It follows easily
that
\eq{}{(\widehat\sV'/\widehat\sV)((z')) \xrightarrow{-z'z^2\partial_z+1}
z(\widehat\sW'/\widehat\sW)((z')) 
}
is an isomorphism, proving (i). 
\item[(ii)] The proof here is analogous to Propositions \ref{prop3.7} and
\ref{prop3.9}. The identity $zw=(1-z'z^2\partial_z)(\sum_{r\ge
-N}v_rz'{}^r)$ with $v_r \in \widehat\sV$ forces $v_r=0,\ r<0$, $v_0 =
zw$, and
$v_r = z^2\partial_z v_{r-1}$. Since the slopes are all $>1$, this is a
contradiction unless $w=0$. Thus $M \inj \sF(\infty,\infty)(M)$. Also
$(z^2\partial_z)^{-1}$ is defined and $z$-adically contracting on $M$.
Since this operator intertwines $z'$ on $\sF(\infty,\infty)(M)$, it is
clear that the image of $M$ is stable under taking Laurent series in
$z'$, from which (ii) follows. 
\item[(iii)] and (iv) are straightforward from
(ii). 
\item[(v)] We can assume $\widehat\sV$ and $\widehat\sW$ come by
completion at $z=0$ from a global good lattice pair $\sV, \sW$, and that
these lattices have no higher cohomology, so
\ml{3.15}{\sF(\sM)\otimes k((z')) \cong
\\
\text{coker}\Big(\Gamma(\P^1,
\sV)\otimes k((z')) \xrightarrow{-z'z^2\partial_z+1}
\Gamma(\P^1,z\sW)\otimes k((z'))\Big). 
}
Indeed, let $j:\G_m \inj \P^1$. Consider the diagram
\eq{3.16}{\begin{CD} 0 @>>> \sV @>>> j_*\sM @>>> j_*\sM/\sV @>>> 0 \\
@. @V z^2\partial_z VV @V z^2\partial_z VV @V z^2\partial_z VV \\
 0 @>>> z\sW @>>> j_*\sM @>>> j_*\sM/z\sW @>>> 0. 
\end{CD}
}
By definition of good lattice pair, the arrow on the right is an
isomorphism. Let $p:j_*\sM/\sV \surj j_*\sM/z\sW$ be the natural
surjection. Then $z^2\partial_z - \frac{1}{z'}p: j_*\sM/\sV((z')) \to
j_*\sM/z\sW((z'))$ is easily checked to be an isomorphism. The assertion
in \eqref{3.15} follows by tensoring  \eqref{3.16} with $k((z'))$,
replacing $z^2\partial_z$ with  $z^2\partial_z - \frac{1}{z'}$ in \eqref{3.16} 
and
taking $R\Gamma$, using vanishing for $H^1$ on the left. 

It follows from \eqref{3.13a} and \eqref{3.15} that there is a natural
map $\sF(\sM)\otimes k((z')) \to \sF(\infty, \infty)(M)$. To see
injectivity, an identity of the form
\eq{}{\sum_{i\ge m}w_i z'{}^i = (1-z'z^2\partial_z)\sum_{j\ge n} \hat
v_jz'{}^j 
}
with $w_i \in \Gamma(\P^1, z\sW)$, $\hat v_j \in \widehat\sV$, 
$w_m\neq 0, \hat{v}_n\neq 0$ yields
that $n=m$, and $w_m = \hat v_m \in
\Gamma(\P^1, z\sW)\cap \widehat\sV = \Gamma(\P^1, \sV)$. Then, recursively
$\hat v_j = w_j + z^2\partial_z \hat v_{j-1} \in \widehat\sV \cap
\Gamma(\P^1,z\sW) = \Gamma(\P^1, \sV)$. 

For surjectivity, we assume moreover that $\sW$ is so positive that
$\Gamma(\P^1, z\sW) \surj z\widehat\sW/\widehat\sV$. 
Given $\hat w \in z\widehat\sW$, we
can then find $w\in \Gamma(\P^1, z\sW)$ with $\hat v:= \hat w - w \in
\widehat\sV$. Then
\ml{}{(\hat w - w)z'{}^N = (1-z'z^2\partial_z)\hat vz'{}^N  +
z^2\partial_z\hat vz'{}^{N+1} \in \\
(1-z'z^2\partial_z)\widehat\sV\otimes z'{}^Nk[[z']] + z\widehat\sW
z'{}^{N+1}.  }
Iterating in this fashion, we get a convergent series in $z'$. 
\end{enumerate}
\end{proof}
\begin{rmk}
One can extend Definition \ref{defn3.11} to the case where the slopes 
of $M$ on $k((z))$ 
are $\le 1$. However, when $M$ has slopes $1$, by Lemma \ref{lem3.3} one can take
$\widehat\sV=z\widehat\sW$ and $\sF(\infty, \infty)(M) = 0$. 
\end{rmk}
\begin{prop} \label{prop3.13} With notations and assumptions as above (in particular,
$M^\nabla = (0)$ and the appropriate slope conditions are assumed to
hold for $M$, cf. Propositions \ref{propdefn3.5}, \ref{prop3.7},
\ref{prop3.12})), we have the following rank and irregularity relations for
the local Fourier transforms:
\begin{enumerate}
\item[(i)] $\rm{irreg.}\sF(0,\infty)(M) = \rm{irreg.}(M)$; $\rm{rk}\,
\sF(0,\infty)(M) = \rm{rk}(M) + \rm{irreg.}(M)$.
\item[(ii)]$\rm{irreg.}\sF(\infty,0)(M) = \rm{irreg.}(M)$; $\rm{rk}\,
\sF(\infty,0)(M) = \rm{rk}(M) - \rm{irreg.}(M)$.
\item[(iii)] $\rm{irreg.}\sF(\infty,\infty)(M) = \rm{irreg.}(M)$;
$\rm{rk}\,
\sF(\infty,\infty)(M) = -\rm{rk}(M) + \rm{irreg.}(M)$.
\end{enumerate}
\end{prop}
\begin{proof}\begin{enumerate}
\item[(i)] Let $\widehat\sV, \widehat\sW$ be a good lattice pair for the
connection $M$ over $k((t))$. We have (using
\cite{Mal2}, Thm. 2.1(b) and the properties of good lattices from section
\ref{sec:lft})
$\partial_t:\widehat\sV \cong \frac{1}{t}\widehat\sW$. In particular,
$\partial_t^{-1}$ stabilizes $\widehat\sV$. Under the identification
$\iota: M \cong \sF(0,\infty)(M)$ (cf. Proposition \ref{prop3.7}),
$\partial_t$ corresponds to multiplication by $z'{}^{-1}$, so
$\iota\widehat\sV \subset
\sF(0,\infty)(M)$ is a $k[[z']]$-submodule. Also,
$z'\partial_{z'}\iota\widehat\sV = \iota \partial_t t\widehat\sV = \iota
(t\partial_t +1)\widehat\sV$. Replacing $\widehat\sV$ by $t^N\widehat\sV$
for $N>>0$, we may assume finally that 
\eq{3.19}{\widehat\sW = (t\partial_t+1)\widehat\sV.
}
Indeed, the formal connection $M$ splits,
$M=M_r\oplus M_i$, into regular singular and irregular parts (\cite{Mal},
p. 51, Thm. 2.3). We may assume a similar decomposition for the lattice
pair. To pass from the good lattice condition $\widehat\sW = t\partial_t
\widehat\sV$ to \eqref{3.19} there is no difficulty in the irregular case
because the slopes are $>0$. In the regular singular case, scaling the
lattices with a large power of $t$ eliminates the eigenvalue $-1$, so
\eqref{3.19} holds in that case as well. We have, therefore
$z'\partial_{z'}\iota\widehat\sV = \iota\widehat\sW$, so
$\iota\widehat\sV, \iota\widehat\sW$ are a good lattice pair for
$\sF(0,\infty)(M)$. In particular
\eq{}{\rm{irreg.}(\sF(0,\infty)(M)) = \dim_k
\iota\widehat\sW/\iota\widehat\sV =  \dim_k
\widehat\sW/\widehat\sV = \rm{irreg.}(M). \notag
}

From Proposition-Definition \ref{propdefn3.5}, we have
$\rm{rk}\sF(0,\infty)(M)$ equals the generic rank of $\sF(\sM)$ where
$\sM$ is the Katz extension coinciding with $M$ at $0$. A standard index
calculation, using that the Fourier sheaf has irregularity $1$ at
$\infty$, yields $\rm{irreg.}(M) + \rm{rk}(M)$, as claimed.
A self-contained computation in the spirit of the article is to consider 
a finite dimensional $k$-vector space $U$ as in \eqref{2.6}. It has dimension 
equal to $\frac{1}{t}\widehat\sW/\widehat\sV$, which is the rank of 
$\sF(0,\infty)(M)$, that is 
${\rm dim} \frac{1}{z}\widehat\sW/\widehat \sW + {\rm dim}\widehat\sW/
\widehat\sV={\rm rk}(M) + {\rm irreg.}(M).$
\item[(ii)] Both assertions follow from (i) together with Proposition
\ref{prop3.10}. The direct computation as above shows again that $\iota 
\widehat\sV, \iota \widehat\sW$ is a good lattice pair if $\widehat\sV, \widehat\sW$ is. 
\item[(iii)] Here $M$ is a $k((z))$-connection, and we have a $k$-vector
space isomorphism $\iota: M \cong \sF(\infty,\infty)(M)$. We have
\eq{}{z'\iota \widehat\sW = \iota(z^2\partial_z)^{-1}\widehat\sW \subset
\iota \widehat\sW 
}
so $\iota \widehat\sW \subset \sF(\infty,\infty)(M)$ is a
$k[[z']]$-lattice. Also $z\partial_z:\widehat\sV \cong \widehat\sW$
implies
\ga{}{\iota\widehat\sV = \iota(z\partial_z)^{-1}\widehat\sW =
\iota(z^{-1}z^2\partial_z)^{-1}\widehat\sW =
z'(z'{}^2\partial_{z'})^{-1}\iota \widehat\sW \\
\iota\widehat\sW = (z'\partial_{z'}-1)\iota\widehat\sV. \notag
}
Since in the case of $\sF(\infty,\infty)$ the slopes are assumed $>1$, it
follows that $\iota\widehat\sW = z'\partial_{z'}\iota\widehat\sV$ so
$\iota\widehat\sV, \iota\widehat\sW$ are a good lattice pair for
$\sF(\infty, \infty)(M)$. It follows as in (i) that $M$ and $\sF(\infty,
\infty)(M)$ have the same irregularity. 

Finally, we compute the rank directly. Consider
\eq{}{\widehat\sV[[z']] \xrightarrow{z'\partial_z - 1/z^2}
\frac{1}{z}\widehat\sW[[z']] 
} 
Note $z^{-1}\widehat\sW/z^{-2}\widehat\sV$ is a $k$-vector space of rank
${\rm dim}\ z^{-1}\widehat\sW/z^{-1}\widehat\sV -{\rm dim
}\ z^{-2}\widehat\sV/z^{-1}\widehat
\sV= \text{irreg.}(M) - \text{rk}(M)$. Write $z^{-1}\widehat\sW = S\oplus
z^{-2}\widehat\sV$ for a $k$-vector space $S$ of this rank. It is straightforward to
check that $\sF(\infty,\infty)(M) \cong S((z'))$ as a
$k((z'))$-module, so the rank is the same. 
\end{enumerate}
\end{proof}

\section{Application to rigidity} \label{sec:rigidity} 
The aim of this section is to  apply 
the equivalence of categories proven in section \ref{sec:lft} to the computation of the index
of rigidity of the Fourier transform of holonomic $\sD$ module on
$\P^1$. 

Let us recall the notion of rigidity as defined by N. Katz
(\cite{Ka}, Introduction). Let $U\subset \P^1$ be a non-empty open
set, defined over a field $k$. If $k$ is a finite field and $M$ is a
$\ell$-adic representation on $U$, then $M$ is said to be rigid if it
is uniquely recognized by the induced local $\ell$-adic
representations at the punctures $\P^1\setminus U$. If $k=\C$ and $M$ is
a local system on $U$, then $M$ is said to be rigid if it is uniquely
recognized by its local monodromies at the punctures. The main theorem
proven by N. Katz in \cite{Ka}, Theorem 5.2.1 is that a rigid
 local system or a tame rigid $\ell$-adic representation
 is always obtained from a
rank 1 one after taking convolution and Fourier transform. 
An important technical tool  to prove this fundamental 
classification theorem is the notion of index of irregularity
and the fact that it is 
preserved by Fourier transforms \cite{Ka}, Theorem 3.0.2.
With this, he is then able in the tame case to inductively lower the
rank of the representation by a suitable rank one twist and
convolution. 

Our aim is to show that our construction of local Fourier transforms
and the accompanying equivalences of categories implies
invariance of the index of rigidity by Fourier transform in the 
$\sD$-module case.

Let $X$ be a smooth, complete curve. Recall (\cite{KaExp}, p.
65, Prop. (2.9.8))  that if $M$ is a connection (i.e. a smooth
holonomic $\sD$-module) on an open $\ell:
U\inj X$, then its middle extension $\ell_{!*}M$ sits in an exact
sequence of $\sD$ modules on $X$ 
\ga{4.1}{0\to \ell_{!*}M\to \ell_*M\to \oplus_{x\in X\setminus U}
  (i_x)_*[((M^\vee\otimes \widehat K_x)^\nabla)^\vee]\to 0}
where $i_x: \{x\}\to X$ is the closed embedding,  $i_{x*}$ is the
$\sD$-module direct image, and $\widehat K_x$ is the complete 
local field at $x$.  

\begin{rmk}\label{rmk4.1} The cokernel in \eqref{4.1} can be written 
$$\oplus_x (i_x)_*
H^1_{DR}(\Spec \widehat K_x, M\otimes \widehat K_x).$$
In particular, taking cohomology yields
\ml{}{0 \to H^1(X,\ell_{!*}M) \to  H^1(X,\ell_{*}M) \to \\
\oplus_{x\in 
X\setminus U} H^1_{DR}(\Spec \widehat K_x, M\otimes \widehat
K_x) \to 0. \notag
}
\end{rmk}

 Let $j_\eta: {\rm Spec}(k(X))\to X, j_{\eta, U}: 
{\rm Spec}(k(X))\to U$ be the
inclusions of the generic point.  Then, since $M$ is assumed smooth on
$U$,  
\ga{4.2}{(j_\eta)_{!*} j_{\eta,U}^* M=\ell_{!*}M.}
One defines (see \cite{Ka}, (3.0.2))
\begin{defn} \label{defn4.1}
Let $M$ be a smooth connection on $U\stackrel{j}{\inj} \A^1$. 
Set $\ell=k\circ j$ where $\A^1\stackrel{k}{\inj} \P^1$. 
Then the index of rigidity of $M$ is defined by
${\rm rig}(M)=\chi(\P^1, \ell_{!*} \sE nd(M)).$
\end{defn}
\begin{thm}[Compare \cite{Ka}, Theorem 3.0.3.] \label{thm4.2} 
Let $M$ be a holonomic $\sD$-module on $\A^1$. We assume 
that $M$ as well as its Fourier transform $\sF(M)$ are 
 the
middle 
extensions of their restriction to the generic point. 
Then one has $${\rm rig}(M)={\rm rig}(\sF(M)).$$
\end{thm}
\begin{proof} Once one has established the equivalences 
of category in section \ref{sec:lft}, 
the proof is exactly the same as Katz's proof in the $\ell$-adic case. 
Let us just give an outline. We assume $M$ is smooth on $U
\stackrel{j}{\inj} \A^1$ as above, and we write $T = \A^1\setminus U$. 
\begin{lem}\label{lem4.3} We have an exact sequence ($z_s$ is a local coordinate
at $s$)
\ml{4.3a}{0 \to j_{!*}\sE nd(M) \to  j_{!*}M \otimes j_{!*}M^\vee\to \\
\oplus_{s\in T}  i_{s*}\Big[\sE nd\big(M\otimes k((z_s))\big)^\nabla\Big/ 
\big(M\otimes
k((z_s))\big)^\nabla \otimes  \big(M^\vee\otimes
k((z_s))\big)^\nabla\Big] \to 0
} 
\end{lem}
\begin{proof} The issue is local around the singular points, so we may consider $M$ a
connection on $k((t)),\ j: \Spec k((t)) \inj \Spec k[[t]]$. We have by
(\cite{KaExp} Proposition 2.9.8 p. 65) ($\delta := k((t))/k[[t]] = \sD/\sD t$, where $\sD$
denotes differential operators on $k[[t]]$. The superscript ${}^\vee$ 
means dual in the appropriate sense.)
\eq{4.4}{0 \to j_{!*}M \to j_*M \to ((M^\vee)^\nabla)^\vee\otimes \delta \to
0. 
}
Now $j_{!*}M$ is $\sO$-torsion-free and hence flat. Also $j_*M \otimes
j_{!*}M^\vee = j_*\sE nd(M)$. Replacing $M$ with $M^\vee$ in \eqref{4.4} and
tensoring the resulting sequence with $((M^\vee)^\nabla)^\vee\otimes
\delta$ yields 
\ml{4.5}{((M^\vee)^\nabla)^\vee\otimes 
\delta\otimes j_{!*}M^\vee \cong ((M^\vee)^\nabla)^\vee\otimes_k
(M^\nabla)^\vee \otimes_k\text{Tor}_1^\sO(\delta,\delta) \cong \\
((M^\vee)^\nabla)^\vee\otimes_k
(M^\nabla)^\vee \otimes_k \delta
}
We get, therefore, an exact sequence
\eq{4.6}{0 \to j_{!*}M \otimes j_{!*}M^\vee \to j_* \sE nd(M)  \to 
((M^\vee)^\nabla)^\vee\otimes_k (M^\nabla)^\vee \otimes_k \delta \to 0
}
Now $j_{!*}\sE nd(M)$ is characterized as a sub-$\sD$-module of $j_*\sE nd(M)$
extending $\sE nd(M)$ and having $\text{Hom}_\sD(j_{!*}\sE nd(M), \delta) =
(0)$ (cf. op. cit. Lemma 2.9.1, p. 57). This implies by \eqref{4.6} that
$j_{!*}\sE nd(M) \subset j_{!*}M \otimes j_{!*}M^\vee $. We get a diagram
\minCDarrowwidth.1cm
\eq{4.7}{\begin{CD} 0 @>>> j_{!*}\sE nd(M) @>>> j_{*}\sE nd(M) @>>> \sE nd(M)^\nabla
\otimes \delta @>>> 0 \\
@. @VV\text{injective} V @| @VV \text{surjective} V \\
0 @>>> j_{!*}M \otimes j_{!*}M^\vee @>>> j_{*}\sE nd(M) @>>>
((M^\vee)^\nabla)^\vee\otimes_k (M^\nabla)^\vee \otimes_k \delta @>>> 0
\end{CD}
}
Finally, this yields
\ga{4.8}{0 \to j_{!*}\sE nd(M) \to j_{!*}M \otimes j_{!*}M^\vee \to V\otimes_k
\delta \to 0 \\
V:= \Big[\sE nd(M)^\nabla\Big/(M^\vee)^\nabla\otimes_k M^\nabla\Big]^\vee,
\notag 
}
proving the lemma.\end{proof}

Returning to the proof of the theorem, we deduce from the lemma
\ml{}{{\rm rig}(M)= h^0(\sE nd(M \otimes_{k[t]}k((z))  )) + \chi
(\A^1, j_{!*}\sE nd(F)) = \\
\chi(\A^1, j_{!*}M\otimes j_{!*}M^\vee)
+ h^0(\sE nd(M \otimes_{k[t]}k((z))  )) +\\
\sum_{x\in \A^1\setminus U} h^0(\sE nd(M\otimes_{k[t]}k((t-x)))^\nabla)\Big/
[M\otimes_{k[t]} k((t-x))]^\nabla \otimes [M^\vee \otimes_{k[t]}
k((t-x))]^\nabla.
}
Here we write $t$ for the coordinate on $\A^1$ and $z=\frac{1}{t}$. 

Writing $D$ for Verdier dual, so for example $DM = M^\vee$, we have 
that $D$ commutes with middle extension (op. cit. Corollary 2.9.1.2)
so 
\eq{}{\chi(\A^1, j_{!*}M\otimes j_{!*}M^\vee) = \chi(\A^1,
j_{!*}M\otimes Dj_{!*}M) }
\begin{lem}[cf. \cite{Ka}, Thm. 3.0.4]\label{lem4.4} Let $N$ be a
holonomic
$\sD$-module on
$\A^1$. Then
\eq{}{\chi(\A^1, N\otimes DN) = \chi(\A^1,\sF(N)\otimes D\sF(N))
}
\end{lem}
\begin{proof}Let $\sM$ be a holonomic $\sD$-module on $\A^1$ and let 
$k: \A^1 \inj \P^1$. Define $\pi: k_!\sM \to k_*\sM$ to be the natural
map. By Kashiwara's theorem, $\ker \pi$ and $\text{coker}\ \pi$ have
the form $\delta^{\oplus n}$. It is not hard to show, e.g. by using
the Levelt classification for formal $\sD$-modules, that the $n$ is
the same for $\ker$ and $\text{coker}$, so in particular
\eq{}{\chi(\A^1, \sM) = \chi(\P^1, k_*\sM) = \chi(\P^1, k_!\sM) =:
\chi_c(\A^1, \sM).  
}
Let $D_{-} = [x\mapsto -x]^*\circ D$. Then $D_- \circ \sF = \sF \circ
D$. Also, write $M_1 *_{!+}M_2 := p_!(M_1\boxtimes M_2)$, where $p:
\A^1\times\A^1 \to \A^1$ is the addition map. One has 
\eq{}{\sF(M_1 *_{!+}M_2) = \sF(M_1)\otimes \sF(M_2)
}
(cf.\cite{KaExp}, 12.2.3(5). To get a formula involving $*_{!*}$
one must modify the argument given there, replacing the lower star
pushforward with lower shriek, and the upper shriek pullback with
upper star. By standard theory, the Fourier transform can be computed
either with lower star or lower shriek.) Finally using that lower
shriek commutes with passage to the fibres, we find
\ml{}{\chi_c\Big(\A^1,\sF(N)\otimes D\sF(N)\Big) =
\text{rank}_0\Big(\sF(N)*_{!+} D_-\sF(N)\Big) = \\
\text{rank}_0\Big(\sF(N)*_{!+} \sF(DN)\Big) =
\text{rank}_0(\sF(N\otimes DN)) = \chi_c(N\otimes DN).
}
\end{proof}
 
For a holonomic $\sD$-module $M$ on $\A^1$ with coordinate $t=1/z$,
it will be convenient to write $\widetilde M :=M\otimes k((z))$. As in Katz, we may
write
$\widetilde M = \widetilde M^{\le 1} \oplus \widetilde M^{> 1}$
according to slopes (\cite{Mal}, Thm. 1.5(2), p. 45). One has, since horizontal
endomorphisms respect slopes
\eq{}{h^0\Big(\sE nd(\widetilde M )\Big) = h^0\Big(\sE
nd(\widetilde M^{\le 1})\Big) +
h^0\Big(\sE nd(\widetilde M^{>1})\Big).  
}
One has by Proposition \ref{prop3.12} (iv), (v)
\eq{}{h^0\Big(\sE nd(\widetilde M^{>1})\Big) = 
h^0\Big(\sE nd(\widetilde{\sF M}^{>1})\Big). 
}
One has, by Definition \ref{defn4.1} (here we write $M = j_{!*}F$ for $F$ a
cnnection on $U$)
\ga{4.17}{{\rm rig}(M)= h^0(\sE nd(\widetilde M)) + \chi
(\A^1, j_{!*}\sE nd(F)).}
This yields
\ga{4.18}{{\rm rig}(M)=\chi(\A^1, j_{!*}F\otimes j_{!*}F^\vee)
+ h^0(\sE nd(\widetilde M))+\\ 
+\sum_{x\in \A^1}\dim\Big[ \sE nd_\nabla(\widetilde M)\Big/
\widetilde M^\nabla \otimes (\widetilde M^\vee)^\nabla \Big]
 .\notag}
Since $D$ commutes with middle extensions, we have
\ga{4.19}{\chi(\A^1, M\otimes
D(M)) = \chi(\A^1, j_{!*}F\otimes j_{!*}F^\vee)
} 
where $D(M)$ is the Verdier dual of $M$. By lemma \ref{lem4.4}
\ga{4.20}{  
\chi\Big(\A^1, M\otimes D(M)\Big)=\chi\Big(\A^1, \sF(M)\otimes
D\sF(M)\Big). 
}
We claim the following identities \eqref{4.21} - \eqref{4.23}
\ga{4.21}{
h^0\sE nd\Big(\widetilde M^{>1}\Big)=
h^0 \sE nd\Big(\widetilde{\sF(M)}^{>1}\Big)
}
\ga{4.22}{
h^0\sE nd\Big(\widetilde M^{\le 1}\Big)=\\
\sum_{x\in \A^1}\dim \Big[ \sE nd_\nabla\Big(\widetilde{\sF(M)}\Big)\Big/
\widetilde{\sF(M)}^\nabla \otimes (\widetilde{\sF(M)}^\vee)^\nabla \Big]\notag
}
\ga{4.23}{h^0\sE nd\Big(\widetilde{\sF(M)}^{\le 1}\Big) =\\ 
\sum_{x\in \A^1}
\dim\Big[\sE nd_\nabla\Big(\widetilde M\Big)\Big/ \widetilde M^\nabla
\otimes (\widetilde M^\vee)^\nabla \Big]
.
\notag}
(Here $^{\le 1}$ and $^{>1}$ refer to the slope decomposition.)  
Now \eqref{4.21} follows immediately from 
proposition \ref{prop3.12} iv), v), and \eqref{4.22} is equivalent to
\eqref{4.23}. Finally, \eqref{4.23} follows from Proposition
\ref{prop3.10}, from Remark
\ref{rmk2.6} which implies that for $x\neq x' \in \A^1$, $${\rm
Hom}(\sF(M\otimes_{k[t]} k((t_x))), \sF(M\otimes_{k[t]} k((t_{x'}))))=0,$$
 and from the following lemma.  
\end{proof}
\begin{lem}[Compare  \cite{Ka}, Proposition 3.1.8]\label{lem4.6} Let
  $M$ be a connection on $k((t))$.
Then one has  
\ga{4.10}{\sE nd_\nabla(M)\Big/M^\nabla\otimes (M^\vee)^\nabla = \sE
nd_\nabla\Big(M/M^\nabla \otimes_k \sO\Big).
}
\end{lem}
 \begin{proof}
We consider the isotypical decomposition $\oplus_N M_N$ of $M$, with
${\rm Hom}(N, N')=0$ if $N\neq N'$. 
Let us write it as $M'\oplus M_{\sO}$ with $M'=\oplus M_N$ where this
sum is over $N\neq \sO$. 
Then $h^0(M')=0, {\rm Hom}(M', \sM_{\sO})=0$ thus
 the left hand side (LHS) of \eqref{4.10} fulfills
\ga{4.11}{LHS(M)=LHS(M')+ LHS(M_{\sO}).} 
And the same holds true for the right hand side (RHS)
\ga{4.12}{RHS(M)=RHS(M')+RHS(M_{\sO}).}
Moreover, $LHS(M')=RHS(M')$ as $M'$ has no flat sections. Thus we
reduce the computation to $M=M_{\sO}$, that is $M$ is nilpotent. 
In this case, this is a purely linear algebra problem. We write
$M=\oplus_i M_i$
where $M_i$ is a maximal Jordan block. Then $M_i^\nabla= k$ and we set
$N_i=M_i/\sO$.  
Then \eqref{4.10} is equivalent to saying that for $i,j$, 
one has an exact
sequence
\ga{4.13}{0\to \sO \to {\rm Hom}(M_i, M_j) \to {\rm Hom}(N_i, N_j)\to 0.}

\end{proof}

We recall, in the context of $\sD$-modules, the central beautiful
observation of Katz.
\begin{thm}\label{thm4.6} Let $X$ be a smooth, complete curve. Let $U
\stackrel{j}{\inj} X$ be Zariski open, and let $M$ be an irreducible
connection on $U$. Suppose 
$$\chi(X, j_{!*}\sE nd(M)) \ge 2. $$
Let $M'$ be another irreducible connection on $U$, and assume for all
$x\in X \setminus U$ we have $M\otimes \widehat K_x \cong M'\otimes \widehat K_x$,
where
$\widehat K_x$ is the Laurent series field at $x$. Then $M\cong M'$.
\end{thm}
\begin{proof} The point is that $\chi(X,j_{!*}N)$ for $N$ a connection on
$U$ depends only on $X$ and the $N\otimes \widehat K_x$. (See e.g.
\cite{Mal}, thm. 4.9, p. 69.) In particular, 
\ml{4.28}{2\le \chi(X, j_{!*}\sE nd(M)) = \chi(X, j_{!*}\sH om(M,M')) \\
\le h^0(X, j_{!*}\sH om(M,M')) +h^2(X, j_{!*}\sH om(M,M')). 
}
Since Verdier duality for holonomic $\sD$-modules on a complete smooth
variety commutes with the de Rham functor (cf. \cite{Bor} (5), p. 326), we 
obtain
\eq{4.29}{h^2(X, j_{!*}\sH om(M,M')) = h^0(X,j_{!*}\sH om(M',M)).
}
It follows that at least one of the modules $\sH om(M,M'),\ \sH om(M',M)$ has a
nontrivial horizontal section. By irreducibility, the two modules are
necessarily isomorphic. 
\end{proof}
Recall Katz' definition (\cite{Ka}, Introduction).
\begin{defn} Let $j: U \inj X$ be as above, and let $M$ be an irreducible
connection on $U$. We say that $M$ is rigid if $M'$ an irreducible
connection on $U$ and $M'\otimes \widehat K_x \cong M\otimes \widehat
K_x$ for all $x\in X\setminus U$ implies $M\cong M'$. 
\end{defn}
\begin{cor}[of Theorem \ref{thm4.2}] Let $M$ be a rank $1$ meromorphic
connection on $\P^1$ and assume the slope of $M$ at $\infty$ is $>1$.
Then $\sF(M)$ is rigid.
\end{cor}
\begin{proof}$\sF(M)$ is smooth on $\A^1$ (\cite{Mal}, (1.4)(b), p. 78).
Since $End(M) = \sO$, the hypotheses of Theorem \ref{thm4.2} are
satisfied, and we conclude 
$$\rm{rig}(\sF(M)) = 2. $$
The result now follows from Theorem \ref{thm4.6}. 
\end{proof}

\begin{thm} \label{thm4.10}
Let $X$ be a smooth, complete curve, and let $j:U\inj X$ be a
non-empty open affine. Let $M$ be an irreducible, rigid connection on
$U$. Then $\rm{rig}(M)=2$. 
\end{thm}
\begin{proof}By irreducibility and duality, 
$$h^0(j_{!*}\sE nd(M)) =
h^2(j_{!*}\sE nd(M)) = 1,$$
so the assertion is equivalent to $h^1(j_{!*}\sE nd(M)) = 0$. Let $\sC$ be the
category of augmented, artinian, local $k$-algebras.  Consider the functor $F: \sC
\to {sets}$,
\ml{}{F(R) = \Big \{(\sM,\nabla_{\sM/R})\text{ lifting } (M,\nabla)\ \Big | \\
(\widehat\sM,\nabla_{\widehat\sM/R}) \cong (\widehat M,\nabla_{\widehat
M})\otimes_k R \Big \}\Big/ \text{isom.}
 }
Here the $\ \widehat{}\ $ means restriction to $\widehat U$, the product of power
series fields at points of $X\setminus U$. We will show that this functor is effectively
pro-representable and smooth, with tangent space $H^1(j_{!*}\sE nd(M))$.
Using a criterion of Artin, \cite{A}, we will show there exists a pointed
affine scheme
$(S=\Spec A, 0)$, smooth and of finite type over $k$, and a connection $(N,
\nabla_N)$ on
$U\times S$ relative to $S$ such that \begin{enumerate}
\item $(N,\nabla_N)|_0 \cong (M,\nabla_M)$.
\item $(\widehat N,\nabla_{\widehat N})\cong (\widehat M,\nabla_{\widehat
M})\times S$, where $\widehat N$ denotes the restriction of $N$ to the union of
tubes $\Spec A((t_x))$, where $x \in X\setminus U$ and $t_x$ is a local parameter at $x$. 
\item $\widehat{\sO}_{S,0}$ pro-represents the functor $F$, and $(N,\nabla_N)$ is
universal. 
\end{enumerate}

Assuming for a moment that we have $(N, \nabla_N)$ satisfying these conditions,
consider the connection $H:= \sH om(M\otimes_k \sO_S, N)$. By rigidity, for any point
$s\in S$, the connection on $H\otimes_{\sO_S} k(s)$ has a horizontal section. In
particular, this is the case at the generic point, so there will exist an nonempty
open $T\subset S$ and a horizontal isomorphism $M\otimes_k \sO_T\cong  N|_T$. To
prove the theorem, we need to show that $N\otimes_{\sO_S} k[\epsilon]\cong
M\otimes_k k[\epsilon]$ for any $\tau: \Spec k[\epsilon] \inj S$ centered at $0$.
If $0\in T$ this is clear. If not, we choose a smooth curve $C\inj S$ passing
through $0$ and tangent to $\tau$. We can further assume $C \cap T \neq
\emptyset$. Shrinking $S$ to a neighborhood of $0$, we can assume $0 \in C$ is
defined by $f=0$ and $C\cap T = C\setminus \{0\}$. 
Restricting the above horizontal
isomorphism to $C\setminus \{0\}$ and multiplying by 
a power of $f$, we get a horizontal
injection $i:M \otimes_k \sO_C \inj N|_C$. Since $N|_C$ 
is a coherent sheaf on
$U\times C$ we see that $\cap f^nN|_C = (0)$. Scaling $i$ by an 
appropriate power
of $f$ we can therefore suppose that the restriction to 
$0$, $i_0 : M \to N_0\cong
M$ is not zero. But $M$ is assumed irreducible, so this map is necessarily an
isomorphism. It follows that $M\otimes_k \sO_C \cong N|_C$, so in particular,
$\tau^*N \cong M\otimes_k k[\epsilon]$ as desired.

It remains to show the existence of $S$. Consider a diagram in $\sC$
\eq{}{\begin{CD} R'\times_R R'' @>>> R'' \\
@VVV @VVV \\
R' @>>> R.
\end{CD}
}
Note that, by irreducibility, any flat automorphism of a lifting $\sM$ of
$M$ over $R$ is necessarily constant (i.e. in $R^\times$.) Suppose given
$\sM' \in F(R')$ and $\sM'' \in F(R'')$ which agree in $F(R)$, i.e. there
is a flat isomorphism $\sM'\otimes R \cong \sM''\otimes R$. Fix such an
isomorphism $\phi$. By the above, it is unique upto $R^\times$. Consider a
subsheaf $\sN \subset \sM' \times \sM''$, $\sN = \{(m',m'')\ |\
\phi(m'\otimes R) = m''\otimes R\}$. Clearly, $\sN$ is a relative
connection on $U \times \Spec(R'\times_R R'')$ lifting $\sM'$ and
$\sM''$. Assuming one of the rings $R'$ and $R''$ surjects onto $R$, elements in $R^\times$
lift to say $R'$. We can then modify $\phi$ by an automorphism of
$\sM'$. In this way we see that $\sN \in F(R'\times_R R'')$ is independent
of the choice of $\phi$. Schlessinger's criterion \cite{S}
\eq{}{F(R'\times_R R'') \cong F(R') \times_{F(R)} F(R'')
}
is therefore satisfied. The tangent space is easily computed to be
\eq{}{\ker(H^1_{DR}(U, \sE nd(M)) \to H^1_{DR}(\widehat U, \sE nd(M)) \cong
  H^1(X, j_{!*}\sE nd(M)),
}
(Compare Remark \ref{rmk4.1}.) and it follows again by \cite{S} that $F$ is
prorepresentable. 

Similarly, the obstruction to smoothness lies in
$H^2(X,j_{!*}\sE nd(M))$. Again by irreducibility, the trace map
\eq{}{ H^2(X,j_{!*}\sE nd(M)) \to H^2(X,j_{!*}\sO_U) \cong k
}
is an isomorphism. Clearly, this trace carries the obstruction to lifting
the connection to the corresponding obstruction to lifting the determinant
of the connection. But these determinant connections are parametrized by a
smooth groupscheme so the determinant obstruction vanishes. We conclude
that our deformation functor $F$ is smooth. 

To construct our family $\sM_S$ of connections algebraizing the above formal
moduli, we apply Artin's criterion \cite{A}, Thm. 1.6. 
For this, we need
to show our functor $F$ is {\it effective} and {\it of finite
presentation}. Effectivity means that if $F$ is pro-represented by
$\Lambda$, then there exists $\sM_\Lambda \in F(\Lambda)$ restricting to
the representing object in $\varprojlim F(\Lambda/\mathfrak
m^n_\Lambda)$.  Choose a lattice $\sL \subset \widehat M$ which is stable
under the group of horizontal automorphisms of $\widehat M$. (Let $\sL_0
\subset \widehat M$ be any lattice. Let $e_1,\dotsc,e_n$ be a vector
space basis for the ring of horizontal endomorphisms of $\widehat M$.
Then $\sL := \sum e_i\sL_0$ works.) Let $\sM_n \in F(\Lambda/\mathfrak
m^n)$ be such that $\varprojlim \sM_n\in \varprojlim F(\Lambda/\mathfrak
m^n)$ is universal. By assumption there exist horizontal isomorphisms
$\psi_n : \widehat M \otimes \Lambda/\mathfrak m^n \cong \widehat\sM_n$ at
infinity. We may glue $\sM_n$ on $U\times \Spec(\Lambda/\mathfrak m^n)$
to $\psi_n(\sL\otimes_k \Lambda/\mathfrak m^n)$ to get bundles
$\overline\sM_n$ on $X\times \Spec(\Lambda/\mathfrak m^n)$. Since
$\text{End}_\nabla(\widehat M\otimes \Lambda/\mathfrak m^n) =
\text{End}_\nabla(\widehat M)\otimes \Lambda/\mathfrak m^n$, the
automorphism $\psi_n^{-1}\circ (\psi_{n+1}\otimes \Lambda/\mathfrak m^n)$
stabilizes $\sL\otimes \Lambda/\mathfrak m^n$. It follows that the
$\overline\sM_n$ are compatible. By Grothendieck, there exists
$\overline\sM$ on $X\times \Spec\Lambda$ which induces the
$\overline\sM_n$. The connections on the $\overline\sM_n$ correspond to
splittings of the Atiyah sequences, with bounded poles on $X\setminus U$
corresponding to the fact that the connection on $\widehat M$ does not
stabilize $\sL$. Again, these splittings agree, so we get a connection on
$\overline \sM$. To examine the polar behavior of this connection, let
$m\ge 0$ be such that $\nabla(\sL) \subset \sL(m(X \setminus U))\otimes
\Omega^1_X$. Then the connection on $\overline \sM$ has poles of order
$\le m$ on $(X\setminus U)\times\Spec\Lambda$. Both effectivity and finite
presentation follow from this. 

The existence of $S, \sM_S$ satisfying properties 1-3 above follows from
Artin, proving the theorem. 
\end{proof}
\begin{rmks}\label{rmk4.11} Theorems \ref{thm4.6} and \ref{thm4.10}
together give a cohomological criterion for rigidity of connections. 
This criterion is proven by Katz in \cite{Ka}, section 1, for regular
singular connections, using transcendental methods. It is unknown on the
$\ell$-adic side.  Note also one does not assume $X$ to be $\P^1$ in the
proofs, yet we know (
\cite{Ka}, section 1) that rigidity is meaningful only on $\P^1$. Indeed,
over a curve $X$ of genus $>0$, we can deform $M$ by twisting with a
family of global rank $1$ connections. The sheaf $j_{!*}\sE nd(M)$ contains
$\sO_X$ as a direct summand, so $\text{rig}(M)=2$ implies $X=\P^1$. 
\end{rmks}

\newpage

\bibliographystyle{plain}
\renewcommand\refname{References}

\end{document}